\documentclass[10pt,reqno]{amsart}

\newcommand\version{November 10, 2009}

\usepackage[ansinew]{inputenc}

\usepackage{textcomp}

\usepackage{cite}

\usepackage{amsmath,amsfonts,amsthm,amssymb,amsxtra}
\usepackage[usenames,dvipsnames]{color}

\usepackage{fancyhdr}
\usepackage[scanall]{psfrag}
\usepackage{graphicx}
\usepackage{ifthen}
\usepackage{pstricks}

\newtheorem{theorem}{Theorem}[section]

\newtheorem{lemma}[theorem]{Lemma}
\newtheorem{corollary}[theorem]{Corollary}

\theoremstyle{definition}

\newtheorem{definition}[theorem]{Definition}
\newtheorem{example}[theorem]{Example}

\theoremstyle{remark}

\newtheorem{remark}[theorem]{Remark}

\numberwithin{equation}{section}

\newenvironment{proofof}[1]{\noindent{\it Proof of~#1.}\rm}{\hfill $\Box$}

\renewcommand{\epsilon}{\varepsilon}

\newcommand{\N}{\mathbb{N}}

\newcommand{\R}{\mathbb{R}}

\begin{document}

\title[Inequalities for eigenvalues of quantum graphs --- \version]{On semiclassical and universal inequalities for eigenvalues of quantum
graphs}

\author{Semra Demirel and Evans M. Harrell II}
\address{Evans M. Harrell II, School of Mathematics,
Georgia Institute of Technology,Atlanta GA 30332-0160, USA}
\email{harrell@math.gatech.edu}
\address{Semra Demirel, University of Stuttgart, Department of Mathematics, Institute of Analysis, Dynamics and Modeling, Chair of Analysis and Mathematical Physics, Pfaffenwaldring 57, D-70569 Stuttgart}
\email{Semra.Demirel@mathematik.uni-stuttgart.de}

\begin{abstract}

We study the spectra of quantum graphs with the method of trace identities (sum rules), which are used to derive inequalities of Lieb-Thirring,
Payne-P\'olya-Weinberger, and Yang types, among others. We show that the sharp constants of these inequalities and even their forms depend on the topology of the graph.  Conditions are identified under which the sharp constants are the same as for the classical inequalities; in particular, this is true in the case of trees.  We also provide some counterexamples where the classical form of the inequalities is false. 

\end{abstract}

\maketitle

\section{Introduction}\label{intro}
This article is focused on inequalities for the means, moments, and ratios of eigenvalues of quantum graphs.  A quantum graph is a  metric graph with one-dimensional Schr\"odinger operators acting on the edges and appropriate boundary conditions imposed at the vertices and at the finite external ends, if any. Here we shall define the Hamiltonian $H$ on a quantum graph as the minimal (Friedrichs) self-adjoint extension of the quadratic form 
\begin{equation}
\phi \in C_c^{\infty} \mapsto E(\phi) := \int_{\Gamma}{|\phi^{\prime}|^2 ds},
\end{equation}
which leads to vanishing Dirichlet boundary conditions at the ends of exterior edges and to the conditions at each vertex $v_k$ that $\phi$ is continuous and moreover
\begin{equation}\label{KVC}
\sum_j{\frac{\partial \phi}{\partial x_{kj}}(0^+)} = 0,
\end{equation}
where the sum runs over all edges emanating from $v_k$, and $x_{kj}$ designates the distance from $v_k$ along the $j$-th edge. (Edges connecting $v_k$ to itself are accounted twice.) In the literature these vertex conditions are usually known as Kirchhoff or Neumann conditions.  Other vertex conditions are possible, and are amenable to our methods with some complications, but they will not be considered in this article.  For details about the definition of $H$ we refer to \cite{Kuc}.

Quantum mechanics on graphs has a long history in physics and physical chemistry \cite{Pau,RS}, but recent progress in experimental solid state physics has renewed attention on them as idealized models for thin domains. While the problem of quantum systems in high dimensions has to be solved numerically, since quantum graphs are locally one dimensional their spectra can often be determined explicitly. 
A large literature on the subject has arisen, for which we refer to the bibliography given in \cite{BK,EKKST}.

The subject of inequalities for means, moments, and ratios of eigenvlaues is rather well developed for Laplacians on domains and for Schr\"odinger operators, and it is our aim to determine the extent to which analogous theorems apply to quantum graphs.  For example, when there is a potential energy $V(x)$ in appropriate function spaces,
Lieb-Thirring inequalities provide an upper bound for the moments of the negative eigenvalues $E_j(\alpha)$ of the Schr\"odinger operator 
$ H(\alpha)=-\alpha \nabla^2 +V(x)$ in $L^2(\R^d)$, $\alpha>0$, of the form

\noindent
\begin{equation}\label{LT}
\alpha^{d/2} \sum_{E_j(\alpha)<0} (-E_j(\alpha))^{\gamma} \leq L_{\gamma,d} \int\limits_{\R^d} (V_-(x))^{\gamma+d/2}\,dx
\end{equation}
for some constant $L_{\gamma,d} \geq L_{\gamma,d}^{cl}$, where $L_{\gamma,d}^{cl}$, known as the {\it classical constant}, is given by 
\begin{equation}\nonumber
L_{\gamma,d}^{cl}=\dfrac{1}{(4\pi)^{d/2}} \dfrac{\Gamma(\gamma+1)}{\Gamma(\gamma+d/2+1)}.
\end{equation}
It is known that (\ref{LT}) holds true for various ranges of $\gamma \geq 0$ depending on the dimension $d$; see \cite{LT,Cw,Hun,L,Ro,Wei}. In particular, in \cite{LW} Laptev and Weidl proved that $L_{\gamma,d}= L_{\gamma,d}^{cl}$ for all $\gamma \geq 3/2$ and $d\geq1$, and Stubbe \cite{St} has recently given a new proof of sharp Lieb-Thirring inequalities for $\gamma \geq 2$ and $d\geq1$ by showing monotonicity with respect to coupling constants. His proof is based on general trace identities for operators \cite{HS0,HS1} known as sum rules, which will again be used as the foundation of the present article.

When there is no potential energy but instead the Laplacian is given Dirichlet conditions on the boundary of a bounded domain, then the means of the first $n$ eigenvalues  are bounded from below by the Berezin-Li-Yau inequality in terms of the volume of the domain, and in addition there is a large family of universal bounds on the spectrum, dating from the work of Payne, P\'olya, and Weinberger \cite{PPW}, which constrain the spectrum without any reference to properties of the domain. (For a review of the subject, see \cite{As}.) It turns out that there are far-reaching analogies between 
these ``universal'' inequalities for Dirichlet Laplacians and Lieb-Thirring inequalities, which have led to common proofs based on sum rules \cite{HS0,HS1,HaHe1,HaHe2,St,HS2}. More precisely, some sharp Lieb-Thirring inequalities and some universal inequalities of the PPW family can be viewed as corollaries of a ``Yang-type'' inequality like \eqref{Yang} below, which in turn follows from a sum rule identity.

In one dimension a domain is merely an interval and the spectrum of the Dirichlet Laplacian is a familiar elementary calculation, for which the question of universal bounds is trivial and uninteresting.  A quantum graph, however, has a spectrum that responds in complex ways to its connectedness; if the total length is finite and appropriate boundary conditions are imposed at exterior vertices, then the spectrum is discrete, and questions about 
counting functions, moments, etc. and their relation to the topology of the graph become interesting, even in the absence of a potential energy.
Below we shall prove several inequalities for the spectra of finite quantum graphs, with the aid of the same trace identities we use to derive Lieb-Thirring inequalities.

For Lieb-Thirring inequalities on quantum graphs the essential question is whether a form of \eqref{LT} holds with the sharp constant for $d=1$,
or whether the connectedness of the graph can change the state of affairs. In \cite{EFK} T. Ekholm, R. Frank and H. Kova\v r\'{\i}k proved Lieb-Thirring inequalities for Schr\"odinger operators on regular metric trees for any $\gamma \geq 1/2$, but without sharp constants. We shall show below that trees enjoy a Lieb-Thirring inequality with the sharp constant when $\gamma \ge 2$, but that this circumstance depends on the topology of the graph.  

We begin with some simple explicit examples showing that neither the expected Lieb-Thirring inequality nor the analogous universal inequalities for finite quantum graphs without potential hold in complete generality. As it will be convenient to have a uniform way of describing examples, we shall let $x_{ij}$ denote the distance from vertex $v_i$ along the $j$-th edge $\Gamma_j$ emanating from $v_i$.  We note that every edge corresponds
to two distinct coordinates $x_{ij} = L - x_{i^{\prime}j^{\prime}}$ where $L$ is the length of the edge, and that a homoclinic loop from a vertex $v_i$ to itself is accounted as two edges.

For the operator $-\frac{d^2}{dx^2}$ on an interval, with vanishing Dirichlet boundary conditions, the universal inequality of Payne-P\'olya-Weinberger reduces to $E_2/ E_1 \le 5$, and the Ashbaugh-Benguria theorem becomes $E_2/ E_1 \le 4$, both of which are trivial in one dimension.  But for which quantum graphs do these classic inequalities continue to be valid?  We shall show below that 
the classic PPW and related inequalities can be proved for the case of trees, with Dirichlet boundary conditions imposed at all external ends of edges, using the method of sum rules.  The sum-rule proof does not work for every graph, however, so the question naturally arises whether the topology makes a real difference, or whether a better method of proof is required. The following examples show that the failure of the sum-rule proof in the case of multiply connected graphs is not an artifact of the method but due to a true topological effect.

We refer to graphs consisting of a circle attached to a single external edge as ``simple balloon graphs.'' The external edge may either be infinite or
of finite length with a vanishing boundary condition at its exterior end. Consider first the graph $\Gamma:=\Gamma_1 \cup \Gamma_2$, which consists of a loop $\Gamma_1$  to which a finite external interval $\Gamma_2$ is attached at a vertex  $v_1$.   Without loss of generality 
we may fix the length of the loop as $2 \pi$, while the ``string'' will be of length $L$. \\ \\ \\ \\

\begin{figure}[h!]
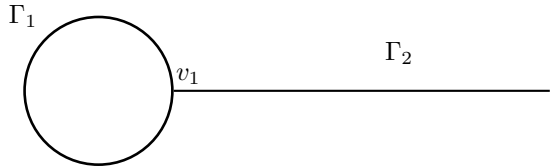

\pscircle[linewidth=1pt](0,0){1}
\qline(6,0)(1,0)
\rput(-1,1){$\Gamma_1$}
\rput(4,0.5){$\Gamma_2$}
\rput(1.2,0.2){$v_1$}
\vspace*{1cm}
\caption{``balloon graph"}
\end{figure}

\begin{example}  
(Violation of the analogue of PPW.)
Let us begin with the case of a balloon graph with $L < \infty$, and assume that there is no potential. We set $\alpha=1$.  Thus $H$ locally has the form $ - \frac{d^2}{dx^2}$ with Dirichlet conditions at the end of the string $\Gamma_2$ and vertex conditions \eqref{KVC} at $v_1$ connecting it to the loop.

For convenience we slightly simplify the coordinate system, letting $x_s := x_{12}$ be the distance on $\Gamma_s=\Gamma_2$ from the node, and $x_{\ell} := x_{11} - \pi$ on $\Gamma_1$.  Thus $x_{\ell}$ increases from $-\pi$ at $v_1$ to $x_2 =+ \pi$ when it joins it again. It is possible to analyze the eigenvalues of the balloon graph quite explicitly: With a Dirichlet condition at $x_s=L,$ any eigenfunction must be of the form
$a \sin(k(L - x_s))$ on $\Gamma_2$.  On $\Gamma_1$ symmetry dictates that the eigenfunction must be proportional to either $\sin{k x_{\ell}}$ or $\cos{k x_{\ell}}$.  There are thus two categories of eigenfunctions and eigenvalues.  Eigenfunctions of the form $\sin{k x_{\ell}}$
contribute nothing to the vertex condition \eqref{KVC} (because the outward derivatives at the node are equal in magnitude with opposite signs), and therefore the derivative of $a \sin(k(L - x_s))$ must vanish at $x_s=0$.  If $k$ is a positive integer, then $k^2$ is an eigenvalue corresponding to an eigenfunction that vanishes on $\Gamma_2$. Otherwise, the conditions on $\Gamma_2$ cannot be achieved without violating the condition of continuity with the eigenfunction on $\Gamma_1$.  To summarize: the eigenvalues of the first category are the squares of positive integers.

The second category of eigenfunctions match $\cos{k x_{\ell}}$ on the loop to $a \sin(k(L - x_s))$ on 
the interval.  The boundary conditions and continuity lead after a standard calculation to the transcendental equation
\begin{equation}\label{transc}
\cot{k L} = 2 \tan{k \pi}.
\end{equation}
There are three interesting situations to consider.  In the limit $L \to 0$, an asymptotic analysis of \eqref{transc} shows that the eigenvalues tend to $\{\left(\frac n 2 \right)^2\}$.  In the limit $L \to \infty$, the lower eigenvalues tend to $\{\left(n + \frac 1 2 \right)^2\frac{\pi^2}{L^2}\}$, which are the eigenvalues of an interval of length $L$ with Dirichlet conditions at $L$ and Neumann conditions at $0$.  The ratio of the first two eigenvalues in this limit is approximately $9$, which is already greater than the classically anticipated value of $5$ or $4$. The highest value of the ratio is, somewhat surprisingly, attained for an intermediate value of $L$, {\it viz}., $L = \pi$, for which \eqref{transc} can be easily solved, yielding $k = \pm \frac 1 \pi \arctan{\frac{1}{\sqrt{2}}} + j$ for a positive integer $j$.  The corresponding fundamental ratio of the lowest two eigenvalues becomes
$$
\dfrac{E_2}{E_1} = \left(\dfrac{\pi - \arctan{\frac{1}{\sqrt{2}}}}{\arctan{\frac{1}{\sqrt{2}}}}\right)^2 \dot{=} 16.8453. 
$$
(We spare the reader the direct calculation showing that the critical value of the ratio occurs precisely at $L = \pi$, establishing this value as the maximum among all simple balloons.)
\end{example}

\begin{example}
(Showing that $E_2 / E_1$ can be arbitrarily large.)  A modification of Example 1.1 with more complex topology shows that no upper bound on the ratio of the first two eigenvalues is possible for the graph analogue of the Dirichlet problem. We again set $\alpha=1$ and assume $V = 0$, and consider a ``fancy balloon'' graph consisting of an external edge, $\Gamma_s$, the ``string,'' of length $\pi$ joined at $v_1$ to $N$ edges $\Gamma_m, m = 1 \dots N$ of length $\pi$, all of which meet at a second vertex $v_2$.  We observe that the eigenfunctions may be chosen either even or odd under pairwise permutation of the edges $\Gamma_m$. This is because if $P f$ represents the linear transformation of a function $f$ defined on the graph by permuting two of the variables $\left\{x_{21}, \dots, x_{2N}\right\}$, and $\phi_j$ is an eigenfunction of the quantum graph with eigenvalue $E_j$, then so are $\phi_j \pm P \phi_j$. (In particular, continuity and \eqref{KVC} are preserved by these superpositions.)
Moreover, the fundamental eigenfunction is even under any permutation, because it is unique and does not change sign.

By continuity and the conditions \eqref{KVC} at the vertices, as in Example 1.1, a straightforward exercise shows that $E_1 = \left(\frac 1 \pi \arctan(\frac{1}{\sqrt{N}})\right)^2$, and that there are other even-parity eigenvalues 
$$\left(j \pm \frac 1 \pi \arctan \left(\frac{1}{\sqrt{N}}\right)\right)^2$$ for all positive integers $j$.  Odd parity, when combined with continuity, forces the eigenfunctions to vanish at the nodes, and thus leads to eigenvalues of the form $j^2$, for positive integers $j$.  The fundamental ratio $E_2 / E_1$ for this example can be seen to be
$$
\left(\frac{\pi - \arctan(\frac{1}{\sqrt{N}})}{\arctan(\frac{1}{\sqrt{N}})}\right)^2,
$$
which is roughly $\pi^2 N$ for large $N$.

\medskip

\noindent
\textbf{Remarks}

1.
With no external edges, the lowest eigenvalue of a quantum graph is $E_1$ = 0, so one might intuitively argue that for a graph with a large and complex interior part the effect of an exterior edge with a boundary condition is small. The theorems and examples given below, however, point towards a more nuanced intuition.

2.
Another instructive example is the ``bunch-of-balloons'' graph, with many nonintersecting loops attached to the string at $v_1$.  We leave the details to the interested reader.
\end{example}

\begin{example}
(Violation of classical Lieb-Thirring.)
Next consider a balloon graph with $L=\infty$ and the Schr\"odinger operator $H(l):=- \frac{d^2}{dx^2}+V(x)$ on $L_2(\Gamma)$ with vertex conditions \eqref{KVC}.
Let the potential $V$ be given by
\[ V(x):=\begin{cases} V_1(x):=\dfrac{-2a^2}{\cosh^2(ax)}\ ,& 
\text{$x\in\Gamma_1=[-\pi,\pi]$}\\[1em] V_2(x):=0\ ,& 
\text{$x\in\Gamma_2=[0,\infty)$} \end{cases} . \]
Then the eigenfunction corresponding to the eigenvalue $-a^2$ is  given by $C \cosh^{-1}(ax_{\ell})$ on $\Gamma_1$ and
by $e^{-ax_s}$ on $\Gamma_2$. The continuity condition gives $C=\cosh(a\pi)$ and the condition \eqref{KVC} at $v_1$ leads to the equation
\begin{equation}\label{bc}
\tanh(a\pi)=\dfrac{1}{2}. 
\end{equation}
Denoting the ratio
$$
\mathcal{Q}(\gamma,V):=\dfrac{|E_1|^{\gamma} } {\int\limits_{\Gamma} |V(x)|^{\gamma+1/2}\,dx},
$$
we compute 
\begin{equation}\nonumber
\mathcal{Q}(3/2,V)= \dfrac{a^3}{2\int\limits_{0}^{\pi} \frac{4a^4}{\cosh^4(ax_{\ell})}\,dx_{\ell}} = \left(8\int\limits_{0}^{a\pi} \frac{1}{\cosh^4(y)}\,dy\right)^{-1} = \left(\frac{8}{3}  \text{tanh}(a\pi)(2+ \text{sech}^2(a\pi)) \right)^{-1}.
\end{equation}
Because of \eqref{bc}, $\text{sech}^2(a\pi) = 1-\text{tanh}^2(a\pi)=\frac{3}{4}$, and therefore
\begin{equation}
\mathcal{Q}(3/2,V)= \dfrac{3}{11} > \dfrac{3}{16}=L_{3/2,1}^{cl}.
\end{equation}
Note that the ratio $\mathcal{Q}(3/2,V)$ is independent of the length of the loop, as expected because any length $L$ can be achieved by a change of scale.\\The ratio $\mathcal{Q}(\gamma,V)$ can also be calculated explicitly for the case $\gamma=2$. In this case
\begin{eqnarray} \nonumber
\mathcal{Q}(2,V)&=&\left[2^{7/2} \left(\dfrac{3}{4} \text{arctan(tanh}(a\pi/2)) + \dfrac{3}{16} \text{sech}(a\pi)+\dfrac{1}{8} \text {sech}^3(a\pi)\right) \right]^{-1}\\ \nonumber
 &\dot{=}& 0.2009 > L_{2,1}^{cl}=\frac{8}{15 \pi} \dot{=}0.1697.
\end{eqnarray}

\end{example}

\section{Lieb-Thirring inequalities for quantum graphs}

\subsection{Classical Lieb-Thirring inequality for metric trees}

Our point of departure is the family of sum-rule identities from \cite{HS0,HS1}. Let $H$ and $G$ be abstract self-adjoint operators satisfying certain mapping conditions. We suppose that $H$ has nonempty discrete spectrum lying below the continuum, $\{E_j:\ H \phi_j=E_j \phi_j\}$.
In the situations of interest in this article the spectrum will either be entirely discrete, in which case we focus on spectral subsets of the form
$J:=\{E_j, j = 1 \dots k\}$, or else, when there is a continuum, it will lie on the positive real axis and we shall take $J$ as the negative part of the spectrum. Let $P_A$ denote the spectral projector associated with  $H$ and a Borel set $A$.

Then, given a pair of self-adjoint operators $H$ and $G$ with domains $D(H)$ and $D(G),$ such that $G(\mathcal J) \subset D(H) \subset D(G)$, where $\mathcal J$ is the subspace spanned by the eigenfunctions $\phi_j$ corresponding to the eigenvalues $E_j$, it is shown in \cite{HS0,HS1} that:
\begin{eqnarray}\label{trace} \nonumber
\quad \sum_{E_j \in J} (z-E_j)^2 \left\langle[G,[H,G]] \phi_j,\phi_j\right\rangle -2(z-E_j)\left\langle[H,G]\phi_j,[H,G]\phi_j\right\rangle \\
=2 \sum_{E_j \in J} \int\limits_{\kappa \in J^c} (z-E_j)(z- \kappa)(\kappa-E_j) \,dG_{j\kappa}^2, 
\end{eqnarray}
where $dG_{j\kappa}^2:=|\left\langle G \phi_j,dP_{\kappa}G \phi_j\right\rangle|$ corresponds to the matrix elements of the operator $G$ with respect to the spectral projections onto $J$ and $J^c$. Because of our choice of $J$,
\begin{eqnarray}\label{HSIneq}
\quad \sum_{E_j \in J} (z-E_j)^2 \left\langle[G,[H,G]] \phi_j,\phi_j\right\rangle -2(z-E_j)
\left\langle[H,G]\phi_j,[H,G]\phi_j\right\rangle \le 0.
\end{eqnarray}
In this section $H$ is the Schr\"odinger operator on the graph $\Gamma$, namely 
\begin{equation}\nonumber
H_{\Gamma}(\alpha)=-\alpha \frac{d^2}{dx^2}+V(x)\ \mbox{in}\ L^2(\Gamma),\ \alpha >0,
\end{equation}
with the usual conditions \eqref{KVC} at each vertex $v_i$. In particular, if any leaves (i.e.\ edges with one free end) are of finite length, vanishing Dirichlet boundary conditions are imposed at their ends. Without loss of generality we may assume that $V\in C_0^{\infty}$ for the operator $H_{\Gamma}(\alpha)$. Under this assumption, for any $\alpha >0$, $H_{\Gamma}(\alpha)$ has at most a finite number of negative eigenvalues. We denote negative eigenvalues of $H_{\Gamma}(\alpha)$ by $E_{ j}(\alpha)$ corresponding to the normalized eigenfunctions $\phi_{j}$.

We shall be able to derive inequalities of the standard one-dimensional type when it is possible to choose $G$ to be multiplication by the arclength along some distinguished subsets of the graph.  This depends on the following:

\begin{lemma}\label{pwiselin}
Suppose that there exists a continuous, piecewise-linear function $G$ on the graph $\Gamma$, such that at each vertex $v_k$
\begin{equation}\label{Gcond}
\sum_j{\frac{\partial G}{\partial x_{kj}}(0^+)} = 0.
\end{equation}
Suppose that $\Gamma = \cup_m \Gamma_m$ with $(G^{\prime})^2 = a_m$ on $\Gamma_m$.  If the spectrum has nonempty essential spectrum,
assume that $z \le \inf \sigma_{\rm ess}(H)$.  Then
\begin{equation}\label{pwiselineq}
\quad \sum_{j,m} (z-E_j)_+^2 a_m \| \chi_{\Gamma_m}\phi_j\|^2 -4 \alpha (z-E_j)_+
a_m  \| \chi_{\Gamma_m}\phi_j^\prime\|^2 \le 0.
\end{equation}

\end{lemma}

\noindent
We observe that $\chi_{\Gamma_m}=1 \Leftrightarrow a_m \neq 0$.

\begin{proof}
The formula \eqref{pwiselineq} is a direct application of \eqref{HSIneq}, when we note that locally, $[H, G] = - 2 G^\prime \frac{d}{dx_{kj}} - G^{\prime \prime}$ and $[G, [H, G]] = 2 (G^\prime)^2$.  (A factor of $2 \alpha$ has been divided out.) The reason for the condition \eqref{Gcond} is that $G \phi_j$ must be in the domain of definition of $H$, which requires that at each vertex,
\begin{equation*}
0 = \sum_j{\frac{\partial G\phi_j}{\partial x_{kj}}(0^+)}
= G \sum_j{\frac{\partial \phi_j}{\partial x_{kj}}(0^+)} + \phi_j \sum_j{\frac{\partial G}{\partial x_{kj}}(0^+)}
= \phi_j \sum_j{\frac{\partial G}{\partial x_{kj}}(0^+)}.
\end{equation*}

\end{proof}

If we are so fortunate that $(G^\prime)^2$ is the same constant on every edge, then \eqref{pwiselineq} reduces to the quadratic inequality
\begin{equation}\label{Yang}
\quad \sum_{j} (z-E_j)_+^2 -4 \alpha (z-E_j)_+  \|\phi_j^\prime\|^2 \le 0,
\end{equation}
familiar from \cite{HS0,HS1,HaHe1,HaHe2,St}, where it was shown that it implies universal spectral bounds for Laplacians and Lieb-Thirring inequalities for Schr\"odinger operators in routine ways. Equation \eqref{Yang} can be considered as a Yang-type inequality, after \cite{Yang}.

\textbf{Stubbe's monotonicity argument}.
In \cite{St} Stubbe showed that some of the classical sharp Lieb-Thirring inequalities follow from the quadratic inequality \eqref{Yang}.  Here we apply the same argument to quantum graphs: 
For any $\alpha >0$, the functions $E_{j}(\alpha)$ are non-positive, continuous and increasing. $E_{j}(\alpha)$ is continuously differentiable except at countably many values where $E_j(\alpha)$ fails to be isolated or enters the continuum. By the Feynman-Hellman theorem,
\begin{equation}\nonumber
\frac{d}{d\alpha} E_{j}(\alpha) = \left\langle\phi_j, - \phi_j^{\prime\prime}\right\rangle
= \|\phi_j^{\prime}\|^2.
\end{equation}
Setting $z=0$, \eqref{Yang} reads
\begin{equation}\nonumber
\alpha \sum_{E_{j}(\alpha)<0}(-E_{j}(\alpha))^2+2 \alpha^2 \frac{d}{d\alpha}\sum_{E_{j}(\alpha)<0}(-E_{j}(\alpha))^2 \leq 0.
\end{equation}
For any $\alpha \in ]\alpha_{N+1},\alpha_N[$ the number of eigenvalues is constant, and therefore
\begin{equation}\nonumber
\frac{d}{d \alpha}\left(\alpha^{1/2}\sum_{E_{j}(\alpha)<0}(-E_{j}(\alpha))^2 \right) \leq 0.
\end{equation}
This means that $\alpha^{1/2}\sum_{E_{j}(\alpha)<0}(-E_{j}(\alpha))^2$ is monotone decreasing in $\alpha$. Hence, by Weyl's asymptotics (see\cite{We1,Bi}),
\begin{equation}\nonumber
\alpha^{1/2} \sum_{E_{j}(\alpha)<0}(-E_{j}(\alpha))^2 \leq \lim_{\alpha \to 0+} \alpha^{1/2} \sum_{E_j (\alpha)<0}(-E_{j}(\alpha))^2 = L_{2,1}^{cl} \int\limits_{\Gamma} (V_-(x))^{2+1/2}\,dx .
\end{equation}
\begin{remark}
Strictly speaking the Feynman-Hellman theorem only holds for nondegenerate eigenvalues. In the case of degenerate eigenvalues one has to take the right basis in the corresponding degeneracy space and to change the numbering if necessary, see e.g. \cite{T}.
\end{remark}

The balloon counterexamples given above might lead one to think that the existence of cycles poses a barrier for a quantum graph to have an inequality of the form \eqref{Yang}.  Consider, however the following example.

\begin{example}  
(Hash graphs.)
Let $\Gamma$ be a planar graph consisting of (or metrically isomorphic to) the union of a closed family of vertical lines and line segments
${\mathfrak F}_v$ and a closed family of horizontal lines and line segments ${\mathfrak F}_h$. We assume that for some $\delta > 0$ the distance between any two lines or line segments in ${\mathfrak F}_v$ is at least $\delta$, and that the same is true of ${\mathfrak F}_h$. (The assumption on the spacing of the lines allows an unproblematic definition of the vertex conditions \eqref{KVC}.) We impose Dirichlet boundary conditions at any ends of finite line segments. We also suppose a ``crossing condition,'' that there are no vertices touching exactly three edges.  (I.e., no line segment from ${\mathfrak F}_v$ has an end point in ${\mathfrak F}_h$ and {\it vice versa}.)  

\noindent
Regarding the graph as a subset of the $xy$-plane, we let $G(x,y) = x+y$.  It is immediate from the crossing condition
that $G$ satisfies \eqref{Gcond}.  Furthermore, the derivative of $G$ along every edge is 1, and therefore the quadratic inequality \eqref{Yang} holds.
\end{example}

A quadratic inequality \eqref{Yang} can arise in a different way, if there is a family of piecewise affine functions $G_{\ell}$ each with a range of values $a_{\ell m}$, but such that $\sum_{\ell}{a_{\ell m}} = 1$ (or any other fixed positive constant).  This occurs in our next example. Even when this is not possible, if we can arrange that $0 < a_{\rm min} \le \sum_{\ell}{a_{\ell m}} \le a_{\rm max}$, then the resulting weaker quadratic inequality
\begin{equation}\label{weak Yang}
\quad \sum_{j} (z-E_j)_+^2 -4 \alpha \frac{a_{\rm max}}{a_{\rm min}}(z-E_j)_+  \|\phi_j^\prime\|^2 \le 0,
\end{equation}
will still lead to universal spectral bounds that may be useful. We speculate about this circumstance below.

\begin{example}($Y$-graph)
As the next example we consider a simple graph, namely the $Y$-graph, which is a 
star-shaped graph with three positive halfaxes $\Gamma_i$, $i=1,2,3$, joined at a single vertex $v_1$.
If we set 
\[ G_1(x):=\begin{cases} g_1:=0\ ,&\text{$x_{11}\in\Gamma_1$}\\[.5em]
 g_2:=-x_{12}\ ,&\text{$x_{12}\in\Gamma_2$}\\[.5em] 
g_3:=x_{13}\ ,&\text{$x_{13}\in\Gamma_3$}\end{cases}  , \]
then obviously $G(\mathcal J) \subset D(H_{\Gamma}(\alpha))$ holds, and with Lemma \ref{pwiselin} we get
\begin{equation}\label{G_1}
\quad \sum_{j} (z-E_j)_+^2 \left( \| \chi_{\Gamma_2}\phi_j\|^2 + \| \chi_{\Gamma_3}\phi_j\|^2 \right) -4 \alpha (z-E_j)_+ \left(\| \chi_{\Gamma_2}\phi_j^\prime\|^2 +  \| \chi_{\Gamma_3}\phi_j^\prime\|^2 \right) \le 0.
\end{equation}
As $\Gamma_1$ doesn't contribute to this inequality, we cyclically permute the zero part of $G$, i.e. we next choose $G_2(x)$, such that $g_2=0,\ g_1=x_{11}$ and $g_3=-x_{13}$, and finally $G_3(x)$, such that $g_3=0,\ g_1=x_{11}$ and $g_2=-x_{12}$. These give us two further inequalities analogous to \eqref{G_1}. Summing all three inequalities, and noting that on every edge, $\sum_{\ell=1}^{3} a_{\ell m} =2$, 
we finally obtain
\begin{equation}
\quad \sum_{j} 2(z-E_j)_+^2 -8 \alpha (z-E_j)_+ \| \phi_j^\prime\|^2 \le 0,
\end{equation}
which when divided by 2 yields the quadratic inequality \eqref{Yang}.
\end{example}

We next extend the averaging argument to prove \eqref{Yang} for arbitrary metric trees.  A metric tree $\Gamma$ consists of a set of vertices, a set of leaves and a set of edges, i.e., segments of the real axis, which connect the vertices, such that there is exactly one path connecting any two vertices.
It is common in graph theory to distinguish between edges and leaves; a leaf is joined to a vertex at only one of its endpoints, ie. there is a free end, at which we shall set Dirichlet boundary conditions. (When the distinction is not material we shall refer to both edges and leaves as edges.  It is also common to regard one free end as the distinguished ``root'' $r$ of the tree, but for our purposes all free ends of the graph have the same status.)
We denote the vertices by $v_i,\ i=1,\dots,n$. The edges including leaves will be denoted by $e$. We shall explicitly write $l_j$ for leaves when the distinction matters.

\begin{theorem}\label{main2}
For any tree graph with a finite number of vertices and edges, the mapping
$$ \alpha \mapsto \alpha^{1/2} \sum_{E_ j(\alpha)<0} (-E_{ j}(\alpha))^2$$
is nonincreasing for all $\alpha > 0$. Consequently
$$ \alpha^{1/2} \sum_{E_{j}(\alpha)<0} (-E_{ j}(\alpha))^2 \leq L_{2,1}^{cl} \int\limits_{\Gamma} (V_-(x))^{2+1/2}\,dx$$
for all $\alpha >0$.
\end{theorem}

\begin{remark}
By the monotonicity principle of Aizenman and Lieb (see \cite{AL}), Theorem \ref{main2} is also true with the sharp constant for higher moments of eigenvalues. Alternatively, the extension to higher values of $\gamma$ can be obtained directly from the trace inequality of \cite{HS2} for power functions with $\gamma > 2$. 
Furthermore, Theorem \ref{main2} can be extended by a density argument to potentials $V \in L^{\gamma+1/2}(\Gamma)$.
\end{remark}

To prepare the proof of Theorem \ref{main2}, we first formulate some auxiliary results.

\begin{lemma}\label{0en=1en}
For all $n \in \mathbb N$,
\begin{equation}\label{0en_auf_e=1en_auf_e}
\sum_{k=0}^{\left[\frac{n-1}{2}\right]} {n-1 \choose 2k} = \sum_{k=0}^{\left[\frac{n}{2}\right]-1}{n-1\choose 2k+1}. 
\end{equation}
 
\end{lemma}
 
\begin{proof}
This is a simple computation.
\end{proof}

\begin{definition}
Let $\mathcal{E}$ be the set of all edges $e \subset \Gamma$. We call the mapping $\mathcal{C}: \mathcal{E} \to \{0,1\}$ a {\it coloring} and say that $\mathcal{C}$ is an {\it admissible coloring} if at each vertex $v\in \Gamma$ the number 
$$\# \{e:\ e\ \mbox{emanates from}\ v:\ \mathcal{C} (e)=1\}$$
is even. We let $\mathcal{A}(\Gamma)$ denote the set of all admissible colorings on $\Gamma$.
\end{definition}

\begin{theorem} \label{coloring}
Let $\Gamma_n$ be a metric tree with $n$ vertices. For an edge $e\subset \Gamma_n$, we denote by
$$a(e,n):=\# \{\mathcal{C}(\Gamma_n) \in \mathcal{A}:\ \mathcal{C}(e)=1\}$$
the number of all admissible mappings $\mathcal{C} \in \mathcal{A}(\Gamma_n)$, such that $\mathcal{C}(e)=1$ for $e\subset \Gamma_n$.
Then
\begin{equation}\label{haupt}
a(e,n)\ \mbox{is independent of}\ e\subset \Gamma_n.
\end{equation}

\end{theorem}

\begin{proof}
We shall prove \eqref{haupt} by induction over the number of vertices of $\Gamma$. The case with one vertex $v_1$ is trivial because of the symmetry of the graph. Given a metric tree $\Gamma_n$ with $n$ vertices, we can decompose it as follows. $\Gamma_n$ consists of a metric tree $\Gamma_{n-1}$ with $n-1$ vertices to which $m-1$ leaves $l_j,\ j=2,\dots,m,$ are attached to the free end of a leaf $l_1 \subset \Gamma_{n-1}$. We call the vertex at which the leaves $l_j,\ j=1,\dots,m,$ are joined $v_n$. Hence,
$$\Gamma_n:= \Gamma_{n-1}\cup v_n \cup \bigcup_{j=2}^{m} l_j.$$ By the induction hypothesis,
\begin{equation}\label{hypo}
a(e,n-1):=
\# \{\mathcal{C} \in \mathcal{A}(\Gamma_{n-1}):\ \mathcal{C}(e)=1\}\ \mbox{is independent of}\ e\subset \Gamma_{n-1}.
\end{equation}
Obviously for every edge or leaf $e \neq l_1$ in $\Gamma_{n-1}$, we have
 \begin{equation}\label{n-1}
a(e,n-1) = \# \{\mathcal{C} \in \mathcal{A}(\Gamma_{n-1}):\ \mathcal{C}(e)=1 \wedge \mathcal{C}(l_1)=1\} +  \# \{\mathcal{C} \in \mathcal{A}(\Gamma_{n-1}):\ \mathcal{C}(e)=1 \wedge \mathcal{C}(l_1)=0\}.
\end{equation}
Now, we have to show that $a(e,n)\ \mbox{is independent of}\ e\subset \Gamma_n$. Note first that for each fixed leaf $l_j$ of the subgraph $\Gamma^*= v_n \cup \bigcup _{j=1}^{m} l_j$, we have
 \begin{equation}
\mu_1:=\# \{\mathcal{C} \in \mathcal{A}(\Gamma^*):\ \mathcal{C}(l_j)=1,\ l_j \in \Gamma^*\} = \sum_{k=0}^{\left[\frac{m}{2}\right]-1}{m-1\choose 2k+1}
\end{equation}
and
\begin{equation}
\mu_0:=\# \{\mathcal{C} \in \mathcal{A}(\Gamma^*):\ \mathcal{C}(l_j)=0,\ l_j \in \Gamma^* \} = \sum_{k=0}^{\left[\frac{m-1}{2}\right]} {m-1 \choose 2k}.
\end{equation}
Hence, for arbitrary neighboring edges $e',\ e'' \subset \Gamma_{n-1}$ the following equality holds,
\begin{align}\label{mualt1} \nonumber
a(e',n)&= \mu_1 \# \{\mathcal{C} \in \mathcal{A}(\Gamma_{n-1}):\ \mathcal{C}(e')=1 \wedge \mathcal{C}(l_1)=1\}\\
&\quad + \mu_0 \# \{\mathcal{C} \in \mathcal{A}(\Gamma_{n-1}):\ \mathcal{C}(e')=1 \wedge \mathcal{C}(l_1)=0\},
\end{align}
and respectively
\begin{align}\label{mualt0} \nonumber
a(e'',n) &= \mu_1 \# \{\mathcal{C} \in \mathcal{A}(\Gamma_{n-1}):\ \mathcal{C}(e'')=1 \wedge \mathcal{C}(l_1)=1\}\\
&\quad + \mu_0 \# \{\mathcal{C} \in \mathcal{A}(\Gamma_{n-1}):\ \mathcal{C}(e'')=1 \wedge \mathcal{C}(l_1)=0\}.
\end{align}
By Lemma \ref{0en=1en}, $\mu:=\mu_0 =\mu_1$. Therefore, with \eqref{n-1} the equalities \eqref{mualt1} and \eqref{mualt0} read
$$
a(e',n)=\mu a(e',n-1),
$$
$$
a(e'',n)=\mu a(e'',n-1).
$$
Furthermore, by the induction hypothesis,
$$
a(e',n-1)=a(e'',n-1),
$$
from which it immediately follows that
\begin{equation*}
a(e',n) = \mu a(e',n-1) = \mu a(e'',n-1) = a(e'',n).
\end{equation*}
This proves Theorem \ref{coloring}.
\end{proof}

\begin{proofof}{Theorem \ref{main2}}
In order to apply Stubbe's monotonicity argument \cite{St}, we need to establish inequality \eqref{Yang} for metric trees. To do this, we proceed as for the example of the $Y$-graph. Let $\mathcal{J}$ denote the subspace spanned by the eigenfunctions $\phi_{j}$ on $L^2(\Gamma)$ corresponding to the eigenvalues $E_{j}$. Note first that there exist self-adjoint operators $G$, which are given by piecewise affine functions $g_i$ on the edges (or leaves) of $\Gamma$, such that $G(\mathcal{J}) \subset D(H(\alpha)) \subset D(G).$  Edges (or leaves) on which constant functions $g_i$ are given, 
do not contribute to the sum rule. Therefore we average over a family of operators $G$, such that every edge $e$ (or leaf) of the tree appears equally often in association with an affine function having $G^{\prime} = \pm 1$ on $e$.
We let $\mathcal{G}$ denote the set of continous operators $G(x)=\{g_i(x)\ \mbox{affine},\ x\in e_i\ (\textit{or}\ l_i)\}$, which satisfy \eqref{KVC} at the vertices $v$ of $\Gamma$. Indeed it is not necessary to average over all the operators $G\in \mathcal{G}$, because it makes no difference in Lemma \ref{pwiselin}, for instance, whether $g'_i=1$ or $g'_i=-1$. Therefore we define an equivalence relation $\sim_G$ on $\mathcal{G}$ as follows: Let $\tilde{G}=\{\tilde{g}_i(x)\ \mbox{affine},\ x\in e_i,\ (\textit{or}\ l_i)\}$ be another operator in $\mathcal{G}$. We say that $G \sim \tilde{G} \Leftrightarrow \forall i \in\{1,\dots,n\}:\ |g'_i(x)|=|\tilde{g}'_i(x)|.$ We define $\mathcal{G}^*:= \mathcal{G}/ \sim$. Then we can consider the isomorphism
\begin{equation}
\mathcal{I}: \mathcal{A}(\Gamma) \to \mathcal{G}^*,
\end{equation}
where for each $ \mathcal{C} \in \mathcal{A}(\Gamma)$ we choose an affine function $G_{\mathcal{C}} \in\mathcal{G}^*$ on $\Gamma$, such that $|G'_{\mathcal{C}} (e)|=\mathcal{C}(e)$ for every $e\subset \Gamma$ . By Theorem \ref{coloring}, we know that $\# \{\mathcal{C} \in \mathcal{A}(\Gamma):\ \mathcal{C}(e)=1\}\ \mbox{is independent of}\ e\subset \Gamma$.  
This means that summing up all inequalities corresponding to  \eqref{pwiselineq}, which we get from each $G_{\mathcal{C}} \in\mathcal{G}^*$, leads to 
\begin{equation}
\quad \sum_{j} (z-E_j)_+^2  p -4 \alpha (z-E_j)_+  p \|\phi_j^\prime\|^2 \le 0,
\end{equation}
\noindent
where $p:=\sum_{\ell} a_{\ell m}=\# \{\mathcal{C} \in \mathcal{A}(\Gamma):\ \mathcal{C}(e)=1\}$ and we have used the normalization $\|\phi_j\|=1$. Having the anologue of inequality \eqref{Yang} for metric trees, we can reformulate the monotonicity argument for our case. This proves Theorem \ref{main2}.
\end{proofof}

\begin{remark}
The proof applies equally to metric trees with leaves of infinite lengths.
\end{remark}

\subsection{Modified Lieb-Thirring inequalities for one-loop graphs}

In this section we consider the graph $\Gamma$ consisting of a circle to which two leaves are attached. It is not hard to see that the construction leading to Lieb-Thirring inequalities with the sharp classical constant fails for one-loop graphs, because no family of auxiliary functions $G_{\ell}$ exists with the side condition that $\sum_{\ell}a_{\ell m}  = 1$ throughout $\Gamma$. Unlike the case of the balloon
graph, it is possible to replace the classical inequality with a weakened version \eqref{weak Yang} as mentioned above.  There is, however another option, based on commutators with exponential functions, following an idea of \cite{HS2}: As usual, we define the one-parameter familiy of 
Schr\"odinger operators
\begin{equation}\nonumber
H(\alpha)=-\alpha \frac{d^2}{dx^2}+V(x),\ \alpha > 0,
\end{equation}
in $L^2(\Gamma)$ with the usual conditions \eqref{KVC} at each vertex $v_i$ of $\Gamma$.
The leaves are denoted by $\Gamma_1:=[0,\infty)$ and $\Gamma_2:=[0,\infty)$, while we write $\Gamma_3$ and $\Gamma_4$ for the semicircles 
with lengths $L$. Let $\phi_j$ be the eigenfunctions of $H(\alpha)$ corresponding to the eigenvalues $E_j(\alpha)$.

\begin{theorem}\label{loop}
Let $q:=2\pi /L$. For all $\alpha>0$ the mapping
\begin{equation}\label{R-mononoticity-periodic}
\alpha\mapsto 
{\alpha}^{\frac{1}{2}}\sum_{E_j(\alpha)<0}
    \left(z-\frac{3}{16} \alpha q^2- E_j\right)_{+}^2
\end{equation}
is nonincreasing. Furthermore, for all $z\in\mathbb{R}$ and all $\alpha>0$ the following sharp Lieb-Thirring inequality holds:
\begin{equation}\label{LTPeriodic}
R_2(z, \alpha)
\le \alpha^{-1/2} L_{2,1}^{cl} \int_{\Gamma}{\left(V({x}) - \left(z +
\frac{3}{16} q^2 \alpha \right)\right)_{-}^{2 + 1/2}d{x}},
\end{equation}
where 
$$
R_2(z, \alpha) := \sum_{E_j(\alpha)<z}
    \left(z - E_j(\alpha) \right)_{+}^2.
$$
\end{theorem}

\begin{remark}
Once again, Theorem \ref{loop} can be extended to potentials $V\in L^{\gamma+1/2}(\Gamma)$ and is true for all $\gamma \geq 2$, either by the monotonicity principle of Aizenman and Lieb \cite{AL} or by the trace formula of \cite{HS2} for $\gamma \geq 2$.
\end{remark}

For the proof of Theorem \ref{loop}, we make use of a theorem of Harrell and Stubbe:

\begin{theorem}[{\cite[Theorem 2.1]{HS2}}]\label{AbstrTrace}
Let $H$ be a self-adjoint operator on $\mathcal{H}$, with a nonempty set $J$ of finitely degenerate eigenvalues lying below the rest of the spectrum $J^c$ and $\{\phi_j\}$ an orthonormal set of eigenfunctions of $H$. Let $G$ be a linear operator with domain $\mathcal{D}_G$ and adjoint $G^*$ defined on $\mathcal{D}_{G^*}$
such that $G(\mathcal{D}_H)\subseteq \mathcal{D}_H\subseteq
\mathcal{D}_G$ and $G^*(\mathcal{D}_H)\subseteq
\mathcal{D}_H\subseteq \mathcal{D}_{G^*}$, respectively. Then
\begin{equation}\label{tf2discrete}
\begin{split}
    &\;\frac1{2}\sum_{E_j\in J}  (z-E_j)^2\,\big(\langle[G^*,[H,G]]\phi_j,\phi_j\rangle+\langle[G,[H,G^*]]\phi_j,\phi_j\rangle\big)\\
    &\quad \le \sum_{E_j\in
    J}(z-E_j)\,\left(\|[H,G]\phi_j\|^2+\|[H,G^*]\phi_j\|^2\right).
    \end{split}
\end{equation}

\end{theorem}
\medskip

\begin{remark}
Strictly speaking, in \cite{HS2}
it was assumed that the spectrum was purely discrete.  However, the extension to the case where continuous spectrum is allowed in $J^c$ follows exactly as in Theorem 2.1 of \cite{HS1}.
\end{remark}

\medskip
\begin{proofof}{Theorem \ref{loop}}
In this case it is not possible to get a quadratic inequality from Lemma \ref{pwiselin} without worsening the constants. This follows from the fact that the conditions $\phi_3(0)=\phi_4(0)$ and $\phi_3(L)=\phi_4(L)$ imply that the piecewise linear function $G$ has to be defined equally on $\Gamma_3$ and  $\Gamma_4$. Consequently, the condition \eqref{KVC} can be satisfied only with different values of $a_m$
as in \eqref{weak Yang}, namely $a_1=a_2=4a_3=4a_4$. Our proof of Theorem \ref{loop} consists of three steps. First we apply Lemma \ref{pwiselin}, after which we apply Theorem \ref{AbstrTrace}. Finally we combine both results and apply the line of argument given in \cite{HS2}.

\noindent
\textit{First step:} Using  Lemma \ref{pwiselin} with the choice,
\[G(x):=\begin{cases} g_1:=-2x_{11}\ ,&\text{$x_{11}\in\Gamma_1$}\\[.5em]
 g_2:=2x_{22}+L\ ,&\text{$x_{22}\in\Gamma_2$}\\[.5em]
 g_3:=x_{13}\ ,&\text{$ x_{13}\in\Gamma_3$}\\[.5em] 
 g_4:=x_{14}\ ,&\text{$ x_{14}\in\Gamma_4$}\end{cases}  , \]
we obtain
\begin{align}\label{in1} \nonumber
&4\left(\sum_{E_j(\alpha) <0} (z-E_j(\alpha))_+ ^2 p_{12}(j)-4 \alpha \sum_{E_j(\alpha) <0} (z-E_j(\alpha))_+  p'_{12}(j) \right) \\
&\qquad +\sum_{E_j(\alpha) <0} (z-E_j(\alpha))_+ ^2 p_{34}(j)- 4 \alpha \sum_{E_j(\alpha) <0} (z-E_j(\alpha))_+ p'_{34}(j) \leq 0,
\end{align}
where $p_{ik}(j):= \| \chi_{\Gamma_i} \phi_j\|^2 +  \| \chi_{\Gamma_k} \phi_j\|^2 $ and $p'_{ik}(j):=\| \chi_{\Gamma_i} \phi'_j\|^2 +  \| \chi_{\Gamma_k} \phi'_j\|^2$.
\\\textit{Second step:} Next, in Theorem \ref{AbstrTrace} we set
\[G(x):=\begin{cases} g_1:=1\ ,&\text{$x_{11}\in\Gamma_1$}\\[.5em] 
g_2:=1\ ,&\text{$ x_{22}\in\Gamma_2$}\\[.5em] 
g_3:=e^{-i2\pi x_{13}/L}\ ,&\text{$x_{13}\in\Gamma_3$}\\[.5em] 
g_4:=e^{i2\pi x_{14}/L}\ ,&\text{$ x_{14}\in\Gamma_4$}\end{cases} 
. \]
It is easy to see that $G\phi_j \in D(H_{\alpha})$. With $q:=2\pi/L$, the first commutators work out to be
$$
[H_j,g_j]=0,\ j=1,2,
$$
$$
[H_3,g_3]=e^{-iqx_{13}} \alpha \left(q^2+2iq d/dx \right),\quad [H_4,g_4]=e^{iqx_{14}} \alpha \left(q^2-2iq d/dx \right);
$$
whereas for the second commutators,
\begin{align}
[g^*_j,[H_j,g_j]]&=[g_j,[H_j,g^*_j]]=0,\quad j=1,2,\\ \nonumber
[g^*_j,[H_j,g_j]]&=[g_j,[H_j,g^*_j]]=2 \alpha q^2,\quad j=3,4.
\end{align}
From inequality \eqref{tf2discrete}, we get
\begin{equation}\label{PerID}
\sum_{E_j(\alpha)\in J}  (z-E_j(\alpha))^2 p_{34}(j)
    \le \alpha \sum_{E_j(\alpha)\in
    J}(z-E_j(\alpha))\,\left(q^2  p_{34}(j) + 4 p'_{34}(j)\right).
\end{equation}
\textit{Third step:} Adding \eqref{in1} and \eqref{PerID} we finally obtain
\begin{equation}\label{summe}
2 \left( R_2  (z,\alpha)+2 \alpha \frac{d}{d \alpha} R_2(z,\alpha)\right) 
\leq \alpha q^2 \frac{3}{2} \sum_{E_j\in J} (z-E_j) p_{34}(j),
\end{equation}
or
\begin{equation}\label{summe2}
2 R_2  (z,\alpha)+4 \alpha \frac{d}{d \alpha} R_2(z,\alpha) - \alpha q^2 \frac{3}{2} R_1 \leq 0,
\end{equation}
which is equivalent to 
\begin{equation}\label{scaleddifferential}
\frac{\partial}{\partial \alpha} \left(\alpha^{1/2} R_2(z, \alpha)\right) \le
\frac{3 q^2}{8} \alpha^{1/2} R_1(z, \alpha).
\end{equation}
Letting $U(z, \alpha) := \alpha^{1/2} R_2(z, \alpha)$, the inequality has the form
\begin{equation}\label{PDI}
\frac{\partial U}{\partial \alpha}  \le \frac{3}{16} q^2 \frac{\partial U}{\partial z}.
\end{equation}
Since the expression in \eqref{LTPeriodic} can be written as $U(z - \frac{3}{16} q^2  \alpha, \alpha)$, an application of the chain rule shows that the monotonicity claimed in  \eqref{LTPeriodic} follows from \eqref{PDI}. (We note that \eqref{PDI} can be solved by changing to characteristic variables
$\xi := \alpha - \frac{16z}{3 q^2} $,
$\eta := \alpha + \frac{16z}{3 q^2}$, in terms of which
\begin{equation}\label{U mono}
\frac{\partial U}{\partial \xi}  \le 0.
\end{equation}
I.e., $U$ decreases as $\xi$ increases while $\eta$ is fixed.)  By shifting the variable in
\eqref{U mono}, we also obtain
\begin{equation}\label{UIneq}
U(z,\alpha)  \le U\left(z+ \frac{3}{16} q^2 (\alpha - \alpha_s),\alpha_s\right)
\end{equation}
for $\alpha \ge \alpha_s$.
By Weyl's asymptotics, for all $\gamma\geq 0$,
\begin{equation}\label{sc-limit}
    \underset{\alpha\rightarrow 0+}{\lim}\alpha^{\frac{d}{2}}\;
    \sum_{E_j(\alpha)<z}(z-E_j(\alpha))^{\gamma}=L_{\gamma,d}^{cl}\int_{\Gamma}{\left(V({x}) -
z\right)_{-}^{\gamma + d/2} d{x}},
\end{equation}
see \cite{We1,Bi}.
Hence, as $\alpha_s \to 0$, the right side of \eqref{UIneq} tends to
$$
L_{2,1}^{cl} \int{\left(V({x}) - \left(z + \frac{3}{16}
q^2 \alpha\right)\right)_-^{2 + 1/2} d{x} },
$$
so the conclusion of Theorem \ref{loop} follows.
\end{proofof}
\\
\begin{remark}
Theorem \ref{loop} can be generalized to one-loop graphs to which $2n,\ n\in \N$ equidistant 
halfaxes
are attached.  
\end{remark}

To summarize, in this section we have seen that for some classes of quantum graphs a quadratic inequality \eqref{Yang} can be proved with the classical constants, and that for some other classes of graphs similar statements can be proved at the price of worse constants as in \eqref{weak Yang}, or of a shift in the zero-point energy as in \eqref{LTPeriodic}. 

It is reasonable to ask whether one can look at the connectness of a graph and say whether a weak Yang-type inequality  \eqref{weak Yang} can be proved. As we have seen, this is the case if there exists a family of continuous functions $G_{\ell}$ on the graph such that
\begin{itemize}
\item
On each edge, all the derivatives $\{G_{\ell}^{\prime}\}$ are constant.
\item
At each vertex $v_k$, each function  $G_{\ell}$ satisfies
$$
\sum_j{\frac{d G_{\ell}}{d x_{kj}}(0^+)} = 0.
$$
\item
For each edge $e$ there exists at least one function $G_{\ell}$
with $G_{\ell}^{\prime} \ne 0$.
\end{itemize}

Interestingly, the question of the existence of such a family of functions can be rephrased in terms of the theory of electrical resistive circuits, a subject dating from the mid nineteenth century \cite{Ki}.  We first note that for a suitable family of functions to exist, there must be at least two leaves, which can be regarded as external leads of an electric circuit, bearing some resistance.  (In the finite case let the resistance be equivalent to the length of the leaf, and in the infinite case let it be some fixed finite value, at least as large as the length of any finite leaf.)  Each internal edge is regarded as a wire bearing a resistance equal to the length of the edge.  If we regard the value of $G_{\ell}^{\prime}$ as a current, then Kirchhoff's condition at the vertex of an electric circuit is exactly the condition \eqref{KVC} that $\sum_j{\frac{d G_{\ell}}{d x_{kj}}(0^+)} = 0$, and
the condition that the electric potential $G_{\ell}$ must be uniquely defined at all vertices is equivalent to global continuity of $G_{\ell}$.  It has been known since Weyl \cite{We2} that the currents and potentials in an electric circuit are uniquely determined by the voltages applied at the leads.  There are, however, circuits such that no matter what voltages are applied to the external leads, there will be an internal wire where no 
current flows; the most well-known of these is the Wheatstone bridge.  (See, for instance, the Wikipedia article on the Wheatstone bridge.)

Let us call a metric graph {\it a generalized Wheatstone bridge} when the corresponding circuit has exactly two external leads and a configuration for which no current will flow in at least one of its wires. Then we conjecture that there are only two impediments to the existence of a suitable family of functions $G_{\ell}$, and therefore to a weakened quadratic inequality \eqref{weak Yang}, namely: Unless a quantum graph contains either
\begin{itemize}
\item{a)} a subgraph that can be disconnected from all leaves by the removal of one point (such as a balloon graph or a graph shaped like the letter $\alpha$); or
\item{b)} a subgraph that when disconnected from the graph by cutting two edges is a generalized Wheatstone bridge,
\end{itemize}
then an inequality of the form \eqref{weak Yang} holds. Otherwise the best that can be obtained may be a modified quadratic inequality with a variable shift, as in Theorem \ref{loop}.
\vspace*{2cm}
\begin{figure}[h!]
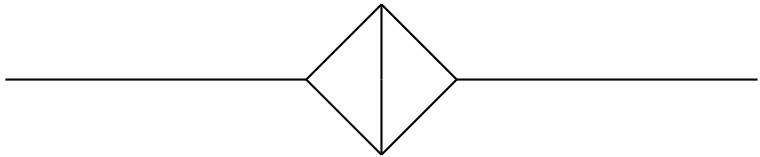

\qline(5,0)(1,0)
\qline(-5,0)(-1,0)
\qline(0,1)(-1,0)
\qline(0,-1)(-1,0)
\qline(0,1)(1,0)
\qline(0,-1)(1,0)
\qline(0,-1)(0,0)
\qline(0,1)(0,0)
\vspace*{1cm}
\caption{``Wheatstone bridge"}
\end{figure}

\section{Universal bounds for finite quantum graphs}

In this section we derive differential inequalities for Riesz means of eigenvalues of the Dirichlet Laplacian on bounded metric trees $\Gamma$
with at least one leaf (free edge). From these inequalities we derive Weyl-type bounds on the averages of the eigenvalues of the Dirichlet Laplacian 
$$
H_D:=\left(-\dfrac{d^2}{dx^2}\right)_D\ \mbox{in}\ L^2(\Gamma),
$$ 
with the conditions \eqref{KVC} at each vertex $v_i$. At the ends of the leaves, vanishing Dirichlet boundary conditions are imposed. We recall that with the methods of \cite{HS0,HaHe1} these are consequences of the same quadratic inequality \eqref{Yang} as was used above to prove Lieb-Thirring inequalities.
When the total length of the graph is finite, the operator $H_D$ on $ D(H_D)$ has a positive discrete spectrum $\{E_j\}_{j=1}^{\infty}$, allowing us to define the {\it Riesz mean of order $\rho$},
\begin{equation} \label{Riesz}
R_{\rho}(z) := \sum_j{(z - E_j)_+^{\rho}}
\end{equation}
for $\rho > 0$ and real $z$.

\begin {theorem} \label{BasicIneq}
Let $\Gamma$ be a metric tree of finite length and with finitely many edges and vertices, and let $H_D$ be the Dirichlet Laplacian in $L^2(\Gamma)$ with domain $D(H_D)$. Then for $z >0$,
\begin{equation} \label{Difference1}
 R_{1}(z) \ge \frac{5}{4z} R_{2}(z);
\end{equation}
\begin{equation} \label{Differential1}
 R_{2}^{\prime}(z) \ge \frac{5}{2z} R_{2}(z);
 \end{equation}
 and consequently
\[\dfrac{R_{2}(z)}{z^{5/2}}\]
is a nondecreasing function of $z$.
\end{theorem}

\begin{proof}
The claims are vacuous for $z \leq E_1$, so we henceforth assume $z>E_1$. The line of reasoning of the proof of Theorem \ref{main2} applies just as well to the operator $H_D$ on $ D(H_D)$, yielding
\begin{equation}\label{YangD}
\quad \sum_{j} (z-E_j)_+^2 -4 (z-E_j)_+  \|\phi_j^\prime\|^2 \le 0.
\end{equation}
Since $V\equiv 0$, $\|\phi_j^\prime\|^2 = E_j$. Observing that
$$
\sum_{j} (z-E_j)_+ E_j = z R_1(z)-R_2(z),
$$
we get from \eqref{YangD}
$$
5 R_2(z) -4z R_1(z) \leq 0.
$$
This proves \eqref{Difference1}. Inequality \eqref{Differential1} follows from \eqref{Difference1}, as $R'_2(z)= 2R_1(z)$.
\end{proof}
Since by the Theorem \ref{BasicIneq}, $R_2(z) z^{-5/2}$ is a nondecreasing function, we obtain a lower bound of the form $R_2(z) \ge C z^{5/2}$ for all $z \ge z_0$ in terms of $R_2(z_0)$. Upper bounds can be obtained from the limiting behavior of $R_2(z)$ as $z \to \infty$, as given by the Weyl law.
In the following, we to follow \cite{HaHe1} to derive Weyl-type bounds on the averages of the eigenvalues of $H_D$ in $L^2(\Gamma)$.

\begin{corollary}\label{twoside}
 For $z\geq 5E_1$,
$$
16 E_1^{-1/2} \left(\frac{z}{5} \right) ^{5/2} \leq R_2(z) \leq L_{2,1}^{cl} |\Gamma| z^{5/2},
$$
where $L_{2,1}^{cl}:= \dfrac {\Gamma(3)}{(4 \pi)^{1/2} \Gamma(7/2)}$, and $|\Gamma|$ is the total length of the tree.
\end{corollary}

\begin{proof}
By Theorem \ref{BasicIneq}, for all $z \geq z_0$,
\begin{equation}\label{mon}
 \dfrac{R_2(z)}{z^{5/2}} \geq  \dfrac{R_2(z_0)}{z_0^{5/2}}.
\end{equation}
As $R_2(z_0) \geq (z_0-E_1)_+^2$ for any $z_0 > E_1$, it follows from \eqref{mon} that
$$
R_2(z) \geq (z_0-E_1)_+^2 \left(\dfrac{z}{z_0}\right)^{5/2}.
$$
The coefficient $\dfrac{(z_0-E_1)_+^2}{z_0^{5/2}}$ is maximized when $z_0=5E_1$. Thus we get
$$
16E_1^{-1/2} \left(\frac{z}{5} \right) ^{5/2} \leq R_2(z).
$$
For metric trees with total length $|\Gamma|$, the Weyl law states that
\begin{equation}\label{weyl}
 \lim_{n\to \infty} \dfrac{\sqrt{E_n}}{n} = \dfrac{\pi}{|\Gamma|},
\end{equation}
(see \cite{Kur}).  It follows that
$$
\dfrac{R_2(z)}{z^{5/2}} \to L_{2,1}^{cl} |\Gamma|,
$$
as $z \to \infty$. Since $\dfrac{R_2(z)}{z^{5/2}}$ is nondecreasing, we get
$$
\dfrac{R_2(z)}{z^{5/2}} \leq  L_{2,1}^{cl} |\Gamma|,\quad \forall z< \infty.
$$
\end{proof}
\noindent
In summary, we get from Theorem \ref{BasicIneq} and Corollary \ref{twoside} the following two-sided estimate:
\begin{equation}\label{summery}
4 E_1^{-1/2} \left(\dfrac{z}{5} \right)^{3/2} \leq \dfrac{5}{4z} R_2(z) \leq R_1(z).
\end{equation}
\medspace
\indent
In order to obtain similar estimates, related to higher eigenvalues, we introduce the notation
$$
\overline{E}_j:=\dfrac{1}{j} \sum_{\ell \le j} E_{\ell}
$$
for the means of eigenvalues $E_{\ell}$; similarly, the means of the squared eigenvalues are denoted
$$
\overline{E^2_j}:=\dfrac{1}{j} \sum_{\ell \le j} E^2_{\ell}.
$$
For a given $z$, we let $ind(z)$ be the greatest integer $i$ such that $E_i \le z$. Then obviously,
$$
R_2(z) = ind(z) (z^2-2z\overline{E_{ind(z)}}+\overline{E^2_{ind(z)}}).
$$
As for any integer $j$ and all $z \geq E_j$, $ind(z) \geq j$, we get
$$
R_2(z) \geq \mathcal{D}(z,j):= j (z^2-2z\overline{E_j}+\overline{E^2_j}).
$$
Using Theorem \ref{BasicIneq} for $z \geq z_j \geq E_j$, it follows that
\begin{equation}\label{D1}
R_2(z)\geq  \mathcal{D}(z_j,j) \left(\dfrac{z}{z_j} \right)^{5/2}.
\end{equation}
Furthermore, $\overline{E_j}^2 \leq \overline{E^2_j}$ by the Cauchy-Schwarz inequality, and hence
\begin{equation}\label{D2}
\mathcal{D}(z,j)=j\left((z-\overline{E_j})^2+ \overline{E^2_j}-\overline{E_j}^2 \right) \geq j(z-\overline{E_j})^2.
\end{equation}
This establishes the following
\begin{corollary}\label{R_1R_2}
Suppose that $z \geq 5 \overline{E_j}$. Then
\begin{equation}\label{R_2}
 R_2(z) \geq \dfrac{16jz^{5/2}}{25 (5 \overline{E_j})^{1/2}}
\end{equation}
and, therefore,
\begin{equation}\label{R_1}
R_1(z) \geq \dfrac{4jz^{3/2}}{5 (5 \overline{E_j})^{1/2}}.
\end{equation}

\end{corollary}

\begin{proof}
Combining equations \eqref{D1} and  \eqref{D2}, we get
$$
R_2(z)\geq  j(z_j-\overline{E_j})^2 \left(\dfrac{z}{z_j} \right)^{5/2}.
$$
Inserting $z_j=5 \overline{E_j}$ the first statement follows. (This choice of $z_j$ maximizes the constant appearing in \eqref{R_2}.) The second statement results from substituting the first statement into \eqref{summery}.
\end{proof}

The Legendre transform is an effective tool for converting bounds on $R_{\rho}(z)$ into bounds on the spectrum, as has been realized previously, e.g., in \cite{LW1}. Recall that if $f(z)$ is a convex function on $\R^+$ that is superlinear in $z$ as $z\to +\infty$, its Legendre transform 
$$
\mathcal{L}[f](w):=\sup_z \{wz-f(z)\}
$$
is likewise a superlinear convex function. Moreover, for each $w$, the supremum in this formula is attained at some finite value of $z$. We also note that if $f(z) \geq g(z)$ for all $z$, then $\mathcal{L}[g](w) \leq \mathcal{L}[f](w) $ for all $w$. The Legendre transform of the two sides of inequality 
\eqref{R_1} is a straightforward calculation (e.g., see \cite{HaHe1}). 
The result is
\begin{equation}\label{LTrafo}
(w-[w])E_{[w]+1}+[w] \overline{E_{[w]}} \leq \dfrac{w^3}{j^2} \dfrac{125}{108} \overline{E_j} ,
\end{equation}
for certain values of $w$ and $j$. In Corollary \ref{R_1R_2} it is supposed that $z \geq 5 \overline{E_j}$. Let $z_{\max}$ be the value for which $\mathcal{L}[f](w)=wz_{\max}-f(z_{\max})$, where $f$ is the right side of \eqref{R_1}. Then by an elementary calculation,
$$
w=\dfrac{6j}{5}\left(\dfrac{z_{\max}}{5\overline{E_j}}\right)^{1/2}.
$$
It follows that inequality \eqref{LTrafo} is valid for $w\geq 6j/5$. Meanwhile, for any $w$ we can always find an integer $k$ such that on the left side of \eqref{LTrafo}, $k-1 \leq w < k$. If $k > 6j/5$ and if we let approach $k$ from below, we obtain from  \eqref{LTrafo}
$$
E_k+(k-1)\overline{E_{k-1}} \leq \dfrac{k^3}{j^2}\dfrac{125}{108} \overline{E_j}.
$$
The left side of this equation is the sum of the eigenvalues $E_1$ through $E_k$, so we get the following:

\begin{corollary}\label{means}
For $k \geq  \frac{6}{5} j$, the means of the eigenvalues of the Dirichlet Laplacian on an arbitrary metric tree with finitely many edges and vertices
satisfy a universal Weyl-type bound,
\begin{equation}\label{125}
 \dfrac{\overline{E_k}}{\overline{E_j}} \leq \frac{125}{108} \left(\frac{k}{j} \right)^2.
\end{equation}

\end{corollary}
In \cite{HS2} it was shown that a similar inequality with a different constant can be proved for all $k \geq j$ in the context of the Dirichlet Laplacian on Euclidian domains. The very same argument applies to quantum graphs with $V=0$. With this assumption $\|\phi'_j\|^2=E_j$, so with $\alpha=1$ \eqref{Yang} can be rewritten as a quadratic inequality,
\begin{equation}\label{polynom}
P_j(z):= \sum_{\ell=1}^{j} (z-E_{\ell})(z-5E_{\ell}) \leq 0
\end{equation}
for $z\in[E_j,E_{j+1}]$ (cf. \cite{HS2}, eq. (4.6)). From \eqref{Difference1} and \eqref{mon} for $z\geq z_0\geq E_j$,
\begin{equation}\label{kette}
R_1(z) \geq \frac{5}{4z} R_2(z) \geq \frac{5}{4} z^{3/2} z_0^{-5/2} \sum_{\ell=1}^j (z_0-E_j)^2.
\end{equation}
The derivative of the right side of \eqref{kette} with respect to $z_0$, by a calculation, is a negative quantity times $P_j(z_0)$, and therefore an optimal choice for the value of \eqref{kette} is the root
\begin{equation}\label{z0}
z_0=3 \overline{E_j} + \sqrt{D_j} \leq 5\overline{E_j},
\end{equation}
where $D_j$ is the discriminant of $P_j$. The inequality in \eqref{z0} results from the Cauchy-Schwarz inequality as in \cite{HS0, HS2}. Because $P_j(z_0)=0$,
$$
0=\sum_{\ell=1}^j (z_0-E_{\ell}) (z_0-5E_{\ell})= 5 \sum_{\ell=1}^j (z_0-E_{\ell})^2-4z_0 \sum_{\ell=1}^j (z_0-E_{\ell}),
$$
so \eqref{kette} reads
$$
R_1(z) \geq \left(\dfrac{z}{z_0} \right)^{3/2} \sum_{\ell=1}^j (z_0-E_{\ell}) = \left(\dfrac{z}{z_0} \right)^{3/2} j (z_0-\overline{E_j}).
$$
From the left side of \eqref{z0}, $z_0-\overline{E_j} \geq \frac{2}{3} z_0$, so
\begin{equation}\label{eh1}
R_1(z) \geq \left(\frac{2}{3} j z_0^{-1/2}\right) z^{3/2}.
\end{equation}
The Legendre transform of \eqref{eh1} is
\begin{equation}\label{eh2}
k \overline{E_k} \leq \dfrac{z_0}{3j^2} k^3,
\end{equation}
and a calculation of the maximizing $z$ in the Legendre transform of the right side of \eqref{eh1} shows that \eqref{eh2} is valid for all $k>j$.
In particular, with the inequality on the right side of \eqref{z0}, we have established the following:

\begin{corollary}
For $k \geq j$, the means of the eigenvalues of $H_D$ in $L_2(\Gamma)$ satisfy
\begin{equation}\label{last}
\dfrac{ \overline{E_k} }{\overline{E_j}} \leq \dfrac{5}{3} \left(\dfrac{k}{j}\right)^2 .
\end{equation}
\end{corollary}
\begin{remark}
Relaxing the assumption to $k \geq j$ comes at the price of making the constant on the right side larger. It would be possible to interpolate between \eqref{last} and \eqref{125} for $k \in [j,6j/5]$ with a slightly better inequality.
\end{remark}

{\bf Acknowledgments}.  The authors are grateful to several people for useful 
comments, including 
Rupert L. Frank, Lotfi Hermi, Thomas Morley, Joachim Stubbe, and Timo Weidl, 
and to Michael Music for calculations and insights generated by them.
We also wish to express our appreciation to the Mathematisches Forschungsinstitut Oberwolfach for 
hosting a workshop in February, 2009, where this collaboration began.
\bibliographystyle{plain}
\bibliography{dh.bib}

\begin{thebibliography}{10}

\bibitem{AL}
Michael Aizenman and Elliott~H. Lieb.
\newblock On semiclassical bounds for eigenvalues of {S}chr\"odinger operators.
\newblock {\em Phys. Lett. A}, 66(6):427--429, 1978.

\bibitem{As}
Mark~S. Ashbaugh.
\newblock The universal eigenvalue bounds of {P}ayne-{P}\'{o}lya-{W}einberger,
  {H}ile-{P}rotter, and {H.} {C}. {Y}ang.
\newblock In {\em Spectral and inverse spectral theory (Goa, 2000)}, Proc.
  Indian Acad. Sci. Math. Sci, 112, pages 3--30. Indian Acad. Sci.

\bibitem{BK}
Gregory Berkolaiko, Robert Carlson, Stephen~A. Fulling, and Peter Kuchment,
  editors.
\newblock {\em Quantum {G}raphs and {T}heir {A}pplications}, volume 415 of {\em
  Contemporary Mathematics}, Providence, RI, 2006. American Mathematical
  Society.

\bibitem{Bi}
Michael~Sh. Birman.
\newblock The spectrum of singular boundary problems.
\newblock {\em Amer. Math. Soc. Trans. (2)}.

\bibitem{Cw}
Michael Cwikel.
\newblock Weak type estimates for singular values and the number of bound
  states of {S}chr\"odinger operators.
\newblock {\em Ann. Math. (2)}, 106(1):93--100, 1977.

\bibitem{EFK}
Thomas Ekholm, Rupert~L. Frank, and Hynek Kovar\'{\i}k.
\newblock Eigenvalue estimates for {S}chr\"odinger operators on metric trees.
\newblock {\em arXiv:0710.5500}.

\bibitem{EKKST}
Pavel Exner, Jonathan~P. Keating, Peter Kuchment, Toshikazu Sunada, and
  Alexander Teplyaev, editors.
\newblock {\em Analysis on graphs and its applications}, volume~77 of {\em
  Proceedings of Symposia in Pure Mathematics}.
\newblock American Mathematical Society, Providence, RI, 2008.
\newblock Papers from the program held in Cambridge, January 8--June 29, 2007.

\bibitem{HaHe2}
Evans~M. Harrell, II and Lotfi Hermi.
\newblock On {R}iesz means of eigenvalues.
\newblock {\em arXiv:0712.4088}.

\bibitem{HaHe1}
Evans~M. Harrell, II and Lotfi Hermi.
\newblock Differential inequalities for {R}iesz means and {W}eyl-type bounds
  for eigenvalues.
\newblock {\em J. Funct. Anal.}, 254(12):3173--3191, 2008.

\bibitem{HS2}
Evans~M. Harrell, II and Joachim Stubbe.
\newblock Trace identities for {C}ommutators with {A}pplications to the
  {D}istribution of {E}igenvalues.
\newblock {\em arXiv:0903:0563v1}.

\bibitem{HS1}
Evans~M. Harrell, II and Joachim Stubbe.
\newblock Universal bounds and semiclassical estimates for eigenvalues of
  abstract {S}chr\"odinger operators.
\newblock {\em arXiv:0808.1133}.

\bibitem{HS0}
Evans~M. Harrell, II and Joachim Stubbe.
\newblock On trace identities and universal eigenvalue estimates for some
  partial differential operators.
\newblock {\em Trans. Amer. Math. Soc.}, 349(5):1797--1809, 1997.

\bibitem{Hun}
Dirk Hundertmark.
\newblock Bound state problems in {Q}uantum {M}echanics.
\newblock In {\em Spectral theory and mathematical physics: a Festschrift in
  honor of Barry Simon's 60th birthday}, Proc. Sympos. Pure Math., LXXVI, part
  1, pages 463--496. Amer. Math. Soc., Providence, R.I., 1980.

\bibitem{Ki}
Gustav~R. Kirchhoff.
\newblock {\"U}ber die {A}ufl\"osung der {G}leichungen, auf welche man bei der
  {U}ntersuchung der linearen {V}ertheilung galvanischer {S}tr\"ome gef\"uhrt
  wird.
\newblock {\em Poggendorf's Ann. Phys. Chemie}, 72, 1847.

\bibitem{Kuc}
Peter Kuchment.
\newblock Quantum graphs: an introduction and a brief survey.
\newblock In {\em Analysis on Graphs and its Applications}, Proc. Symp. Pure.
  Math., AMS 2008, pages 291--314. Amer. MAth. Soc., 2008.

\bibitem{Kur}
Pavel Kurasov.
\newblock Schr\"odinger operators on graphs and geometry. {I}. {E}ssentially
  bounded potentials.
\newblock {\em J. Funct. Anal.}, 254(4):934--953, 2008.

\bibitem{LW1}
Ari Laptev and Timo Weidl.
\newblock Recent results on {L}ieb-{T}hirring inequalities.
\newblock In {\em Journ\'ees ``\'{E}quations aux {D}\'eriv\'ees, {P}artielles''
  ({L}a {C}hapelle sur {E}rdre, 2000)}, pages Exp. \ No. \ XX, 14. Univ.
  Nantes, Nantes, 2000.

\bibitem{LW}
Ari Laptev and Timo Weidl.
\newblock Sharp {L}ieb-{T}hirring inequalities in high dimensions.
\newblock {\em Acta Math.}, 184(1):87--111, 2000.

\bibitem{L}
Elliott~H. Lieb.
\newblock The number of bound states of one-body {S}chr\"odinger operators and
  the {W}eyl problem.
\newblock In {\em Geometry of the {L}aplace operator ({P}roc. {S}ympos. {P}ure
  {M}ath., {U}niv. {H}awaii, {H}onolulu, {H}awaii, 1979)}, Proc. Sympos. Pure
  Math., XXXVI, pages 241--252. Amer. Math. Soc., Providence, R.I., 1980.

\bibitem{LT}
Elliott~H. Lieb and Walter Thirring.
\newblock Inequalities for the moments of the eigenvalues of the
  {S}chr\"odinger {H}amiltonian and their relation to {S}obolev inequalities.
\newblock {\em Studies in Mathematical Physics, Princeton University Press,
  Princeton, NJ}.

\bibitem{Pau}
Linus Pauling.
\newblock The diamagnetic anistropy of aromatic molecules.
\newblock {\em J. Chem. Phys.}

\bibitem{PPW}
L.~H. {P}ayne, G.~{P}\'olya, and H.~F. {W}einberger.
\newblock On the ratio of consecutive eigenvalues.
\newblock {\em J. Math. Phys.}

\bibitem{Ro}
Grigorii~V. Rozenblum.
\newblock Distribution of the discrete spectrum of singular differential
  operators.
\newblock {\em Izv. Vys\v s. U\v cebn. Zaved. Matematika}, (1(164)):75--86,
  1976.

\bibitem{RS}
Klaus Ruedenberg and Charles~W. Scherr.
\newblock Free-electron network model for conjugated systems i, theory.
\newblock {\em J. Chem. Phys.}

\bibitem{St}
Joachim Stubbe.
\newblock Universal monotonicity of eigenvalue moments and sharp
  {L}ieb-{T}hirring inequalities.

\bibitem{T}
Walter Thirring.
\newblock A course in mathematical physics: Quantum mechanics of atoms and
  molecules, pages 149-150.

\bibitem{Wei}
Timo Weidl.
\newblock On the {L}ieb-{T}hirring constants {$L\sb {\gamma,1}$} for
  {$\gamma\geq 1/2$}.
\newblock {\em Comm. Math. Phys.}, 178(1):135--146, 1996.

\bibitem{We1}
H.~Weyl.
\newblock Das asymptotische {V}erteilungsgesetz der {E}igenwerte linearer
  partieller {D}ifferentialgleichungen.
\newblock {\em Math. Ann.}

\bibitem{We2}
Hermann Weyl.
\newblock Repartici\'on de corriente en una red conductora.
\newblock {\em Rev. Mat. Hisp.-Amer}, 5(1):153--164, 1923.

\bibitem{Yang}
Hong~Cang Yang.
\newblock Estimates of the difference between consecutive eigenvalues.
\newblock {\em preprint 1995 (revision of {I}nternational {C}entre for
  {T}heoretical {P}hysics preprint {I}{C}/91/60)}, {T}rieste, {A}pril 1991.

\end{thebibliography}
\end{document}